\documentclass{article}
\usepackage{amssymb}
\usepackage{amsfonts}
\usepackage{mathrsfs}

\usepackage{amssymb, latexsym}

\usepackage[pdftex]{color} 

\usepackage{graphicx}

\usepackage{amsmath}

\definecolor{flesh}{rgb}{.90,.60,.50}
\definecolor{night}{rgb}{.05,0,.05}
\definecolor{gold}{rgb}{.80,.75,.10}

\definecolor{wine}{rgb}{1.00,0,.10}
\definecolor{dimmed}{gray}{0.9}

\definecolor{blau1}{rgb}{0,0.3,0.65}
\definecolor{red1}{rgb}{.5,.0,.0}
\definecolor{green1}{rgb}{.0,.5,.0}
\definecolor{water}{rgb}{.0,.9,0.7}
\definecolor{panelbackground}{rgb}{0,0.39,0.61}
\definecolor{blue1}{rgb}{0.367,0.679,0.839}
\definecolor{blue2}{rgb}{0.000,0.391,0.605}
\definecolor{blue3}{rgb}{0.4,0.8,0.9}

\def\en t{{{\rm Z}\mkern-5.5mu{\rm Z}}}

\newtheorem{theorem}{Theorem}[section]

\newtheorem{conjecture}[theorem]{Conjecture}

\newtheorem{corollary}[theorem]{Corollary}

\newtheorem{lemma}[theorem]{Lemma}

\newtheorem{proposition}[theorem]{Proposition}

\newtheorem{remark}[theorem]{Remark}

\begin{document}

\title{\large\bf  On the Law of Large Numbers for the empirical measure process of Generalized Dyson Brownian motion
}

\author{Songzi Li, \ \ Xiang-Dong Li\thanks{Research supported by NSFC No. 11371351, Key Laboratory RCSDS, CAS, No. 2008DP173182, and a
Hundred Talents Project of AMSS, CAS.}, \ \ Yong-Xiao Xie}


\maketitle

\begin{center}
\begin{minipage}{120mm}
\begin{center}{\bf Abstract}\end{center}  We study  the generalized Dyson Brownian motion (GDBM) of
an interacting $N$-particle system with logarithmic Coulomb
interaction and general  potential $V$. Under reasonable condition
on $V$, we prove the existence and uniqueness of strong solution to
SDE for GDBM. We then prove that the family of the empirical
measures of GDBM is tight on $\mathcal
 {C}([0,T],\mathscr{P}(\mathbb{R}))$ and all the large $N$ limits
satisfy a nonlinear McKean-Vlasov equation. Inspired by  previous
works due to Biane and Speicher, Carrillo, McCann and Villani, we
prove that the McKean-Vlasov equation is indeed the gradient flow of
the Voiculescu free entropy on the Wasserstein space of probability
measures over $\mathbb{R}$. Using the optimal transportation theory,
we prove that if $V''\geq K$ for some constant $K\in \mathbb{R}$,
the McKean-Vlasov equation has a unique weak solution. This proves
the Law of Large Numbers and the propagation of chaos for the
empirical measures of GDBM. Finally, we prove the longtime
convergence of the McKean-Vlasov equation for $C^2$-convex
potentials $V$.
\end{minipage}
\end{center}


\section{Introduction}

In 1962, F. Dyson \cite{Dy1, Dy2} observed  that the eigenvalues of
the $N\times N$ Hermitian matrix valued Brownian motion is an
interacting  $N$-particle system with the logarithmic Coulomb
interaction and derived their statistical properties. Since then,
the Dyson Brownian motion has been used in various areas in
mathematics and physics, including statistical physics and the
quantum chaotic systems. See e.g. \cite{Meh} and reference therein.
In \cite{RS93}, Rogers and Shi proved that the empirical measure of
the eigenvalues of the $N\times N$ Hermitian matrix valued
Ornstein-Uhlenbeck process weakly converges to the nonlinear
McKean-Vlasov equation with quadratic external potential as $N$
tends to infinity. This also gave a dynamic proof of Wigner's famous
semi-circle law for Gaussian Unitary Ensemble. See also \cite{AGZ,
Gui}.

The purpose of this paper is to study the generalized Dyson Brownian
motion  and the associated McKean-Vlasov equation with the
logarithmic Coulomb interaction and with general external potential.
More precisely, let $\beta\geq 1$ be a parameter, $V:
\mathbb{R}\rightarrow \mathbb{R}^+$ be a continuous function, let
$(W^{1},\ldots, W^{N})$ be an $N$-dimensional Brownian motion
defined on a filtered probability space~$(\Omega, \mathscr{F},
(\mathscr{F}_t)_{t\geq 0}, \mathbb{P})$ satisfying the usual
conditions. Let
~$\lambda_N(0)=(\lambda^1_N(0),\ldots,\lambda^N_N(0))\in
\bigtriangleup_{N}=\{(x_{i})_{1\leq i\leq N}\in\mathbb{R}^{N}:
x_{1}<x_{2}<\ldots<x_{N}\}$. The generalized Dyson Brownian motion
$({\rm GDBM})_V$ is an interacting $N$-particle system
$\lambda_N(t)=(\lambda^1_N(t), \ldots, \lambda_N^N(t))$  with the
logarithmic Coulomb interaction and with external potential $V$, and
is defined as the solution to the following SDEs
\begin{eqnarray} d\lambda_N^i(t)=\sqrt{\frac{2}{\beta
N}}dW^i_{t}+\frac{1}{N}\sum\limits_{j:j\neq i}
\frac{1}{\lambda^i_N(t)-\lambda^j_N(t)}dt-{1\over 2}V'(\lambda_N^i(t))dt, \ \ \ i=1, \ldots, N, \label{SDE1}
\end{eqnarray}
with initial data ~$\lambda_N(0)$. It is a SDE for $N$-particles
with a singular drift of the form ${1\over x-y}$ due to the
logarithmic Coulomb interaction, and an additional nonlinear drift
due to non quadratic external potential. When $V=0$ and $\beta=1, 2,
4$, it is the standard Dyson Brownian motion \cite{Dy1, Dy2}. When
$V(x)={x^2\over 2}$ and $\beta>1$, it has been studied by Chan
\cite{Ch}, Rogers and Shi \cite{RS93}, C\'epa and L\'epingle
\cite{CL}, Fontbona \cite{JF1, JF2}, Guionnet \cite{Gui}, Anderson,
Guionnet and Zeitouni \cite{AGZ} and references therein. When $N=\infty$, see \cite{KT}. 

By It\^o's calculus, $({\rm GDBM})_V$ is an interacting $N$-particle
system with the Hamiltonian
$$H(x_1, \ldots, x_N):=-\frac{1}{2N}\sum\limits_{1\leq i\neq
j\leq N}\log |x_i-x_j|+\frac{1}{2}\sum\limits_{i=1}^NV(x_i),$$
and the infinitesimal generator of $({\rm GDBM})_V$
 is given by
\begin{eqnarray*}
\mathscr{L}^\beta_{N}f=\frac{1}{\beta
N}\sum\limits_{k=1}^N\frac{\partial^2f}{\partial
x_k^2}+\sum\limits_{k=1}^{N}\left({\rm P.V.}
\int_{\mathbb{R}}\frac{L_N(dy)}{x_k-y}-\frac{1}{2}V'(x_k)\right)\frac{\partial
f}{\partial x_k},
\end{eqnarray*}
where  $f\in C^2(\mathbb{R}^N)$ and
$L_N=\frac{1}{N}\sum\limits_{i=1}^N\delta_{x_i}\in\mathscr{P}(\mathbb{R})$.

Under suitable condition on $V$, we prove that the SDEs
$(\ref{SDE1})$ for $({\rm GDBM})_V$ admit a unique strong solution
$\lambda_N(t)\in \bigtriangleup_{N}$ with infinite lifetime. See
Theorem \ref{Th1} below. Let
$$L_N(t)=\frac{1}{N}\sum\limits^N_{i=1}\delta_{\lambda^i_N(t)}\in\mathscr{P}(\mathbb{R}), \ \ \ t\in [0, \infty).$$Standard argument shows that the family
$\{L_N(t), t\in [0, T]\}$ is tight on $C([0, T],
\mathscr{P}(\mathbb{R}))$, and the limit of any weakly convergent
subsequence of $L_N(t)$, denoted by $\mu_t$, is a weak solution to
the following nonlinear McKean-Vlasov  equation: for all $f\in
C^2_b(\mathbb{R})$,
\begin{eqnarray}\label{DBM7}
\frac{d}{dt}\int_{\mathbb{R}}
f(x)\mu_t(dx)=\frac{1}{2}\int\int_{\mathbb{R}^2}\frac{\partial_xf(x)-\partial_yf(y)}{x-y}\mu_t(dx)\mu_t(dy)-\frac{1}{2}\int_{\mathbb{R}}
V'(x)f'(x)\mu_t(dx).
\end{eqnarray}
In the case $\mu_t$ is absolutely continuous with respect to the
Lebesgue measure on $\mathbb{R}$, integrating by parts, one can
verify that the probability density $\rho_t={d\mu_t\over dx}$
satisfies the following nonlinear McKean-Vlasov equation (also called nonlinear Fokker-Planck equation in the literature) 

\begin{eqnarray}
{\partial \rho_t\over \partial t}={\partial \over \partial
x}\left(\rho_t\left({1\over 2}V'-{\rm
H}\rho_t\right)\right),\label{NFK1}
\end{eqnarray}
where
$$
{\rm H}\rho(x)={\rm P.V.}\int_{\mathbb{R}}{\rho(y)\over x-y}dy$$ is
the Hilbert transform of $\rho$.

It seems that one can not find well-established result in the literature
on the uniqueness of weak solutions to the above McKean-Vlasov
equation with general external potential $V$. By lack of this, one
can not find established result in the literature on the Law of
Large Numbers for the GDBM with non quadratic potentials. One of the
main observations of this paper is to find (and prove) the fact
that the McKean-Vlasov equation is indeed the gradient flow of the
Voiculescu free entropy $\Sigma_V$ on the Wasserstein space
$\mathscr{P}_2(\mathbb{R})$ equipped with Otto's infinite
dimensional Riemannian structure, and to use the optimal
transportation theory to prove the uniqueness of weak solutions to
the McKean-Vlasov equation for general potentials $V$ with natural
condition. This allows us to further derive the Law of Large Numbers
for the empirical measures of
 the generalized Dyson Brownian motion.

 Following Voiculescu
\cite{VD1}and Biane \cite{Bian03}, for every $\mu\in
\mathscr{P}(\mathbb{R})$, we introduce the Voiculescu free entropy
as follows
\begin{eqnarray*}
\Sigma_V(\mu)=-\int_{\mathbb{R}}\int_{\mathbb{R}}
\log|x-y|d\mu(x)d\mu(y)+\int_{\mathbb{R}}V(x)d\mu(x).
\end{eqnarray*}
By \cite{Joh98}, it is well-known that if $V$ satisfies the growth condition
\begin{eqnarray}
V(x)\geq (1+\delta)\log(x^2+1),\ \ \ \ x\in \mathbb{R},\label{grow}
\end{eqnarray}
then there exists a unique minimizer (called the equilibrium
measure) of $\Sigma_V$, denoted by
\begin{eqnarray*}
\mu_V={\rm arg min}_{\mu\in \mathscr{P}(\mathbb{R})} \Sigma_V(\mu).
\end{eqnarray*}
Moreover, $\mu_V$ satisfies the Euler-Lagrange equation
\begin{eqnarray*}
{\rm H}\mu_V(x)={1\over 2}V'(x),\ \ \ \ \forall x\in \mathbb{R}.
\end{eqnarray*}
The relative free entropy is defined as follows
\begin{eqnarray*}
\Sigma_V(\mu|\mu_V)=\Sigma_V(\mu)-\Sigma_V(\mu_V).
\end{eqnarray*}
Following \cite{VD1, Bian03}, the relative free Fisher information
is defined as follows
\begin{eqnarray*}
{\rm I}_V(\mu)=\int_{\mathbb{R}}\left({\rm H}\mu(x)-{1\over
2}V'(x)\right)^2d\mu(x).
\end{eqnarray*}
Note that
\begin{eqnarray*}
{\rm I}_V(\mu_V)=0.
\end{eqnarray*}

We now state the main results of this paper. Our first result
establishes the existence and uniqueness of the strong solution to
SDEs $(\ref{SDE1})$ and the tightness of the associated empirical
measure for a class of $V$ with reasonable condition. 

\begin{theorem}\label{Th1} \footnote{Under the condition $-xV'(x)\leq C$ for all $x\in \mathbb{R}$,  Rogers and Shi \cite {RS93} 
proved the non-collision of the strong solution to $(\ref{SDE1})$, but they did not precisely state the condition $(i)$ which is need for the existence of solution. In \cite{GM},  Graczyk and Malecki proved the existence and uniqueness of strong 
solution to SDE $(\ref{SDE1})$ under the assumption that $V'$ is global Lipschitz. The conditions in Theorem \ref{Th1} require that $V'$ satisfies the local monotonicity condition, i.e., $(i)$,  and one-side growth condition at infinity, i.e., $(ii)$. 
We would like to point out that the  local monotonicity condition $(i)$ in Theorem \ref{Th1} is weaker than the condition $V'$ is local Lipschitz, and the one-side growth condition $(ii)$ is 
also weaker than the condition $V'$ is global Lipschitz . } Let $V$ be a $C^1$ function satisfying the growth condition $(\ref{grow})$ and the following conditions\\
(i) For all $R>0$, there is $K_R>0,$~such that for all
~$x,y\in\mathbb{R}$ with $|x|, |y|\leq R$, $$(x-y)(V'(x)-V'(y))\geq
-K_R|x-y|^2,$$ (ii) There exists a constant $\gamma>0$ such that
\begin{eqnarray}
-xV'(x)\leq \gamma(1+|x|^2), \ \ \forall \ x\in
\mathbb{R}.\label{Cond1}
\end{eqnarray}
Then, for all $\beta\geq 1$, and for any given $\lambda_N(0)\in
\bigtriangleup_{N}$, there exists a unique strong solution
$(\lambda_N(t))_{t\geq 0}$ taking values in $\bigtriangleup_{N}$
with infinite lifetime to SDEs $(\ref{SDE1})$ with initial value
$\lambda_N(0)$.\\

Moreover, suppose that
$L_N(0)\rightarrow\mu\in\mathscr{P}(\mathbb{R})$ as
$N\rightarrow\infty$, and
\begin{eqnarray*}
\sup\limits_{N\geq 0}\int_{\mathbb{R}}\log(x^2+1)dL_N(0)<\infty.
\end{eqnarray*}
Then, the family $\{L_N(t), t\in [0, T]]\}$ is tight in $\mathcal
 {C}([0,T],\mathscr{P}(\mathbb{R}))$, and the limit of any weakly convergent
subsequence of $\{L_N(t), t\in [0, T]]\}$ is a weak solution of the
McKean-Vlasov equation $(\ref{DBM7})$.
\end{theorem}

 Inspired by previous works due to Biane \cite{Bian03},
 Biane-Speicher \cite{BS01},  Carrillo-McCann-Villani \cite{CMV} (see Theorem 3.1 below) and Sturm \cite{Stu}, we can prove the following result which might be already known by experts even though we cannot find the explicit statement in the literature.

\begin{theorem}\label{GFV}\label{MVGF}  For all $V: \mathbb{R}\rightarrow [0,
\infty)$ being a $C^1$ function satisfies the condition
$(\ref{Cond1})$, the nonlinear McKean-Vlasov equation
$(\ref{NFK1})$, i.e.,
\begin{eqnarray*}
{\partial \rho_t\over \partial t}=-{\partial \over \partial
x}(\rho_t({\rm H}\rho_t-\frac{1}{2}V'))
\end{eqnarray*}
is indeed the gradient flow of \ $\Sigma_V$ on the Wasserstein space
$\mathscr{P}_2(\mathbb{R})$.
\end{theorem}

In the optimal transportation theory, it is well known that if the free energy $F$ on the Wasserstein space  is $K$-convex, then the 
$W_2$-Wasserstein distance between the solutions of the gradient flow $\partial_t \mu=-{\rm grad} F(\mu)$ with initial datas $\mu_1(0)$ and $\mu_2(0)$ satisfies $W_2(\mu_1(t), \mu_2(t))\leq e^{-Kt}W_2(\mu_1(0), \mu_2(0))$. 
See \cite{Ot, OV, Stu, StR, Vi1, Vi2}. 
In view of this and Theorem \ref{GFV}, and using the Hessian calculation for nonlinear diffusions with interaction on the Wasserstein space as developed in \cite{CMV, Stu}, we can prove the following result, which ensures the uniqueness
of weak solutions to the McKean-Vlasov equation with general 
potential $V$ satisfying the condition $V''\geq K$.

\begin{theorem}\label{Th2}
 Suppose that $V$ is a $C^2$ function satisfying the same condition as in Theorem \ref{Th1}, and there exists a constant $K\in \mathbb{R}$ such that
$$V''(x)\geq K, \ \ \ \ \forall\ x\in \mathbb{R}.$$
Then the Voiculescu free entropy  $\Sigma_V$ on the Wasserstein space $\mathscr{P}_2(\mathbb{R})$ is $K$-convex, i.e., its Hessian on $\mathscr{P}_2(\mathbb{R})$ satisfies 
\begin{eqnarray*}
{\rm Hess}_{\mathscr{P}_2(\mathbb{R})}\Sigma_V\geq K.
\end{eqnarray*}
Let $\mu_i(t)$ be two solutions of the McKean-Vlasov equation
$(\ref{NFK1})$ with initial data $\mu_i(0)$, $i=1, 2$. Then for all
$t>0$, we have
\begin{eqnarray*}
W_2(\mu_1(t), \mu_2(t))\leq e^{-Kt}W_2(\mu_1(0), \mu_2(0)).
\end{eqnarray*}
In particular, the Cauchy problem of the McKean-Vlasov equation
$(\ref{NFK1})$ has a unique weak solution.
\end{theorem}

We would like to point out that C\'spa and L\'epingle  \cite{CL} proved the uniqueness of weak solution to the McKean-Vlasov equation $(\ref{NFK1})$ with quadratic potential function $V(x)=ax^2+bx$ with two constants, $a\geq 0$ and $b\in \mathbb{R}$,  and 
Fontbana \cite{JF2} proved the uniqueness of weak solution to the McKean-Vlasov equation $(\ref{NFK1})$ with external potential $V$ such that 
$V'(x)=\theta x+b_1(x)$, where $\theta\in \mathbb{R}$ is a constant and $b_1\in C^1(\mathbb{R})$ is a bounded function with bounded derivative. See also \cite{JF1}. Theorem \ref{Th2} establishes the uniqueness of weak solution to the 
McKean-Vlasov equation $(\ref{NFK1})$ with more general external potentials $V$ satisfying the  condition $V''\geq K$  for some $K\in \mathbb{R}$. 

As a consequence of Theorem \ref{Th1} and Theorem \ref{Th2}, we can
derive the Law of Large Numbers for the empirical measures of the
generalized Dyson Brownian motion.

\begin{theorem}\label{Th3}\footnote{For $V(x)=Kx^2$ with $K\in \mathbb{R}$, the result in Theorem \ref{Th3} also holds for $p=2$}
 Suppose that $L_N(0)$ weakly converges to $\mu(0)\in \mathscr{P}(\mathbb{R})$. Let $V$ be a $C^2$ function satisfying the same condition as in Theorem \ref{Th1} and $V''\geq K$ for some constant $K\in \mathbb{R}$. Then the empirical measure $L_N(t)$ weakly converges to the unique solution $\mu_t$ of the McKean-Vlasov equation $(\ref{NFK1})$. Moreover, for all $p\in [1, 2)$, we have
\begin{eqnarray*}
W_p(\mathbb{E}(L_N(t), \mu_t)\rightarrow 0 \ \ \ {\rm as}\ \ N\rightarrow \infty,\label{speed1}
\end{eqnarray*}
where the convergence is uniformly with respect to $t\in [0, T]$ for
all fixed $T>0$.
\end{theorem}

The notion of propagation of chaos, which was introduced by M. Kac,
plays a critical role in the study of the large $N$ limit of
$N$-particle systems. According to Sznitman-Tanaka's theorem
\cite{ST}, for exchangeable systems, propagation of chaos is
equivalent to the law of large numbers for the empirical measures of
the system. In view of this and Theorem \ref{Th3}, we have the
following result, which is a dynamic version of a result due to Johansson (Theorem 2 in \cite{Joh98}). 

\begin{theorem}\label{THPC}
Assume the conditions in Theorem \ref{Th3} holds. Let $M_{N;k}(t;dx_1,\cdots,dx_k)$ be the $k$-th moment measure for the random probability measure $L_N(t,\cdot)$, that is, for any Borel sets $A_1,\cdots,A_k$,
\begin{eqnarray*}
M_{N;k}(t; A_1,\cdots,A_k):=\mathbb{E}(L_N(t,A_1)\cdots L_N(t,A_k)).
\end{eqnarray*}
Then we have
\begin{eqnarray*}
\lim\limits_{N\rightarrow \infty}\int_{\mathbb{R}^k}\varphi(x_1,\cdots,x_k) M_{N;k}(t;dx_1,\cdots,dx_k)=\int_{\mathbb{R}^k}\varphi(x_1,\cdots,x_k)\mu_t(dx_1)\cdots\mu_t(dx_k)
\end{eqnarray*}
for any continuous, bounded $\varphi$ on $\mathbb{R}^k$.
\end{theorem}

By the ergodic theory of SDE, for a wide class of potentials $V$,
and for any fixed $N$, it is known that $L_N(t)$ converge to $L_N$,
as $t\rightarrow \infty$. On the other hand, the large $N$-limit of
$L_N(t)$, i.e., $\mu_t(dx)=\rho_t(x)dx$, satisfies the nonlinear
Fokker-Planck equation $(\ref{NFK1})$. It is natural to ask the
question whether $\mu_t$ converges to $\mu_V$ in the weak
convergence topology or with respect to the $W_2$-Wasserstein
distance for general potentials $V$.  If this is true, then, with
respect to the weak convergence on $\mathscr{P}(\mathbb{R})$ or the
$W_2$-Wasserstein topology on $\mathscr{P}_2(\mathbb{R})$, the
following diagram is commutative
\begin{eqnarray*}
L_N(t)&\Longrightarrow& \mu_t\\
\Downarrow& & {\Downarrow }\\
L_N&\Longrightarrow &\mu_V
\end{eqnarray*}
In other words, we have
\begin{eqnarray*}
\lim\limits_{N\rightarrow \infty}\lim\limits_{t\rightarrow
\infty}L_N(t)= \lim\limits_{t\rightarrow
\infty}\lim\limits_{N\rightarrow \infty}L_N(t).
\end{eqnarray*} In the literature, Chan  \cite{Ch} and Rogers-Shi \cite{RS93} proved that this is true
for $V(x)={x^2\over 2}$. See also \cite{AGZ, Gui}. In particular,
this gives a dynamic proof of Wigner's semi-circle law for the
Gaussian Unitary Ensemble. The following result provides some
positive answers to this problem for $C^2$-convex potentials.

\begin{theorem} \label{Th5}(i) Suppose that $V$ is $C^2$-convex, i.e., $V''\geq 0$. Then $\mu_t$ converges to $\mu_V$ with respect to the Wasserstein distance in $\mathscr{P}_2(\mathbb{R})$, i.e.,
\begin{eqnarray*}
W_2(\mu_t, \mu_V)\rightarrow 0 \ \ \ {\rm as}\ \  t\rightarrow
\infty.
\end{eqnarray*}
(ii) Suppose that $V$ is $C^2$ and there exists a constant $K\in
\mathbb{R}$ such that
$$V''(x)\geq K, \ \ \ \forall x\in \mathbb{R}.$$
Then for all $t>0$, we have
\begin{eqnarray*}
\Sigma_V(\mu_t|\mu_V)&\leq& e^{-2Kt}
\Sigma_V(\mu_0|\mu_V),\\
W_2(\mu_t, \mu_V)&\leq& e^{-Kt}W_2(\mu_0, \mu_V).
\end{eqnarray*}
In particular,  if $V$ is $C^2$-uniform convex with $V''\geq K>0$, then $\mu_t$ converges to $\mu_V$  with the exponential rate $K$ in the $W_2$-Wasserstein topology on $\mathscr{P}_2(\mathbb{R})$. \\
(iii) Suppose that $V$ is a $C^2$, convex and there exists a
constant $r>0$ such that
$$V''(x)\geq K>0, \ \ \ \ |x|\geq r.$$
Then $\mu_t$
converges to $\mu_V$ with an
exponential rate in the $W_2$-Wasserstein topology on  $\mathscr{P}_2(\mathbb{R})$. More precisely, there exist two constants $C_1>0$ and $C_2>0$ such that
\begin{eqnarray*}
W^2_2(\mu_t, \mu_V)\leq {e^{-C_1t}\over C_2}\Sigma_V(\mu_0|\mu_V), \ \ \ \ \ t>0.
\end{eqnarray*}
\end{theorem}

As a corollary of Theorem \ref{Th5}, for $C^2$-convex potentials, we
can give a dynamic proof of the well-known result due to Boutet de
Monvel-Pastur-Shcherbina \cite{BPS} and Johansson \cite{Joh98}.
Their result says that, for $V$ satisfying the growth condition
$(\ref{grow})$, the empirical measure $L_N={1\over
N}\sum\limits_{i=1}^N \delta_{x_i}$ weakly converges to the
equilibrium measure $\mu_V$, where $(x_i, i=1, \ldots, N)$,
satisfies the following probability distribution
\begin{eqnarray*}
P^N_\beta(dx_1,\ldots,dx_N)=\frac{1}{Z^\beta_N} \Pi_{i\neq
j}|x_i-x_j|^{\frac{\beta}{2}}\exp\left(-\frac{\beta N
}{2}\sum\limits_{i=1}^NV(x_i)\right)\prod_{i=1}^N dx_i,
\end{eqnarray*}
where $\beta>0$ is a parameter. We would like to mention that, for
non-convex potentials $V$, we do not know how to give a dynamic
proof of the above result. We would like to mention a recent paper by  Bourgade, Erd\"os, and Yau \cite{BEY}  in which the authors proved
 the bulk universality of the $\beta$-ensembles with non-convex regular analytic potentials for any $\beta>0$. Whether or not their idea of introducing a ``convexified measure" can be used to 
 extend the results in this paper to non-convex case, will be an interesting problem for study in future. 

Finally, let us mention that, for $\beta=2$ and for real analytic
function $V$, we can prove that the generalized Dyson Brownian
motion can be realized as the eigenvalues process of the $N\times N$
real Hermitian matrix valued diffusion process defined by
\begin{eqnarray*}
dX^N_t={1\over \sqrt{N}}dB^N_t-{1\over 2}V'(X^N_t)dt,\label{DBM1}
\end{eqnarray*}
where $B^N_t$ is the $N\times N$ Hermitian matrix valued Brownian
motion. Moreover, we can prove that $X_t^N$ converges in
distribution to the free diffusion process $X_t$, which was defined
by Biane and Speicher \cite{BS01}. This extends a famous result, due
to Voiculescu \cite{VD2, VD3} and Biane \cite{Bian95}, which states
that the renormalized Hermitian Brownian motion ${1\over
\sqrt{N}}B^N_t$ converges in distribution to the free Brownian
motion $S_t$. See \cite{LX2014}.

The rest of this paper is organized as follows. In Section $2$, we
prove Theorem \ref{Th1}. In Section $3$, we prove Theorem \ref{GFV},
Theorem \ref{Th2} and Theorem \ref{Th3}, Theorem \ref{THPC}. In
Section $4$, we prove Theorem \ref{Th5}. In Section $5$, we discuss
the case of double-well potential and raise some conjectures. Finally, let us mention that this paper is an update revised version of our previous paper entitled Generalized Dyson Brownian motion,
McKean-Vlasov equation and  eigenvalues of random matrices (arxiv.org/abs/1303.1240v1).

\section{Proof of Theorem \ref{Th1}}

The proof of Theorem \ref{Th1}  is adapted from classical argument coming back to McKean and exposed in \cite{RS93, CL, AGZ}. \\

\noindent{\bf Proof of existence and uniqueness of GDBM}. First, for
fixed $R>0$, let~$\phi_R(x)=
        x^{-1}$~  if~$|x|\geq R^{-1}$,~
        and~$\phi_R(x)=R^2x$ if $|x|< R^{-1}$. Since ~$\phi_R$~is uniformly Lipschitz and $V$ satisfies $(i)$ and $(ii)$, by Theorem 3.1.1 in \cite{Ro}, the following SDE for the truncated Dyson Brownian
motion
\begin{equation}
 d\lambda^i_{N,R}(t)=\sqrt{\frac{2}{\beta
N}}dW^i_{t}+\frac{1}{N}\sum\limits_{j:j\neq i}
\phi_R(\lambda^i_{N,R}(t)-\lambda^j_{N,R}(t))dt-\frac{1}{2}V'(\lambda^i_{N,R}(t))dt,\label{E2.1}
\end{equation}
with $\lambda^i_{N,R}(0)=\lambda^i_N(0)$~for~$1\leq i\leq N$,
 has a unique strong
solution. Let
\begin{eqnarray*}
\tau_R:=\inf\{t:\min_{i\neq
j}\mid\lambda^i_{N,R}(t)-\lambda^j_{N,R}(t)\mid<R^{-1}\}.
\end{eqnarray*}
Then $\tau_R$ is monotone increasing in $R$ and $\lambda_{N,R}(t)=\lambda_{N,R'}(t)$ for all $t\leq
\tau_R$ and $R<R'$.

Second, let $\lambda_N(t)=\lambda_{N, R}(t)$ on $t\in [0, \tau_R)$.
To prove that $\lambda_N(t)$ is a global solution
 to SDE $(\ref{SDE1})$, we need
only to prove $\lambda_N(t)$ does not explode, and  $\lambda_N^i(t)$ and $\lambda_N^j(t)$ never collide for all $t>0$, $i\neq j$.

To prove that $\lambda_N(t)$ does not explode, let
$R_t=\frac{1}{2N}\sum\limits_{j=1}^N\lambda_N^j(t)^2$. By It\^{o}'s
formula, and by Levy's characterization, we can introduce a new
Brownian motion~$B$,~such that ~$$dR_t={2\over N}\sqrt{R_t\over
\beta}dB_t+\left(\frac{1}{\beta N}+\frac{N-1}{2N}-\frac{1}{2}\langle
L_N(t), xV'(x)\rangle\right)dt.$$ Let $R'$ be the solution of
$$dR'_t={2\over N}\sqrt{R'_t\over \beta}dB_t+\left(\frac{1}{\beta
N}+\frac{N-1}{2N}+\frac{1}{2}\gamma+\gamma R'_t\right)dt,$$ with
$R'_0=R_0.$ Under the assumption $(\ref{Cond1})$, and using the
comparison theorem of one dimensional SDEs, cf. \cite{IW}, we can
derive that
$$R_t\leq
R'_t,\ \ \ \forall \ t\geq 0,~~{\rm a.s}.$$Moreover, by Ikeda and
Watanabe \cite{IW} (p. 235-237), the process $R'$ never explodes. So
the process $R$ (and hence $\lambda_N(t)$) ~does not explode in
finite time \footnote{ In \cite{RS93},  Rogers and Shi proved the
non-explosion of GDBM for $V$ satisfying $-xV'(x)\leq \gamma$,
$\forall x\in \mathbb{R}$. }.

To prove that $\lambda_N^i(t)$ and $\lambda_N^j(t)$ never collide for all $t>0$, $i\neq j$, let us introduce the Lyapunov function $f(x_1,\ldots,x_N)=\frac{1}{N}\sum\limits_{i=1}^N
V(x_i)-\frac{1}{N^2}\sum\limits_{i\neq j}\log|x_i-x_j|$. Similarly to \cite{Gui, AGZ}, we can prove
\begin{eqnarray*}
df(\lambda_N(t))&=&dM_N(t)+\frac{1}{N^3}\left(\frac{1}{\beta}-1\right)\sum\limits_{k\neq
i}\frac{1}{(\lambda_N^i(t)-\lambda_N^k(t))^2}dt-\frac{1}{2N}\sum\limits_{i=1}^N|V'(\lambda_N^i(t))|^2dt\\
& & +\frac{1}{N^2}\left({1\over \beta}\sum\limits_{i=1}^N
V''(\lambda_N^i(t))+\frac{3}{2}\sum\limits_{j\neq
i}\frac{V'(\lambda_N^i(t))-V'(\lambda_N^j(t))}{\lambda_N^i(t)-\lambda_N^j(t)}\right)dt,
\end{eqnarray*}
where~$M_N$~is the following local
martingale~$$dM_N(t)=\frac{2^{\frac{1}{2}}}{\beta^{\frac{1}{2}}
N^{\frac{3}{2}}}\sum\limits_{i=1}^N \left(
V'(\lambda_N^i(t))-\frac{1}{N}\sum\limits_{k:k\neq
i}\frac{1}{\lambda_N^i(t)-\lambda_N^k(t)} \right)dW_t^i.$$
Fix $K>0$ and $R>0$ such that $\lambda_N^i(0)\in
[-K, K]$ and $|\lambda_N^i(0)-\lambda_N^j(0)|\geq R^{-1}$ for all $i\neq j$, $i,
j=1, \ldots, N$.  Let $C_1(K)\geq 0$ be such that $\sup\limits_{x\in [-K, K]}V''(x)\leq
C_1(K)$. Let $A_N(t)dt=df(\lambda_N(t))-dM_N(t)$, and
$\zeta_K=\inf\limits\{t\geq 0: \lambda_N^i(t)\notin [-K, K],\ {\rm
for\ some}\ i=1, \ldots, N\}$, then for any fixed $T>0$, $\sup\limits_{t\in [0, T]}A_N(t\wedge \zeta_K)\leq C_1(K)$ and
$\{f(\lambda_N(t\wedge \zeta_K)-C_1(K)(t\wedge \zeta_K),\ \ t\in [0,
T]\}$ is a supermartingale. Let $C_2(K):= \inf\limits\{V(x):|x|\leq K\}$, we can prove
\begin{eqnarray*}
\mathbb{P}(\tau_R\leq \zeta_K\wedge T)\leq
\frac{N^2(f(\lambda_N(0))+TC_1(K))+N(N-1)\log(2K)-C_2(K)}{\log(2K)+\log R}.
\end{eqnarray*}
Letting $R$, $T$ and $K$ tend to infinity, we can prove $\mathbb{P}(\tau_\infty<\zeta)=0$, where $\zeta:=\inf\{t:\lambda_N^i(t)=\lambda_N^j(t)~\exists~1\leq i\neq j\leq N\}$.
 This proves that $\lambda_N^1(t), \ldots, \lambda_N^N(t)$ does not collide.

Finally, by the continuity of the trajectory of $\lambda_N(t)$, we
have $\lambda_N(t)\in\triangle_N$ for all $t\geq 0$. The same
argument as used in the  proof of Theorem 12.1 in \cite{Gui} proves
the uniqueness of the weak solution to SDEs $(\ref{SDE1})$.  The
proof of Theorem \ref{Th1} is completed.

\medskip

\noindent{\bf Proof of tightness and identification of McKean-Vlasov
limit}

We follow the argument used in \cite{RS93} to prove the tightness of
$\{L_N(t),t\in[0,T]\}$. Let us pick functions $f_j\in
C_b^\infty(\mathbb{R}, \mathbb{C}), j=1,2,\ldots,$ which is dense in
$C_b(\mathbb{R})$. Thus
$$\langle\mu, f_j\rangle=\langle\mu', f_j\rangle, \ \ \forall
j\Rightarrow  \mu=\mu'.$$ We also pick a $C^\infty$ function
$f_0:\mathbb{R}\rightarrow[1, \infty)$ with the properties
$$f_0(x)=f_0(-x),\ \ \ f_0(x)\rightarrow\infty \ \ \ {\rm as} \ \ x\rightarrow\infty, \ x\in \mathbb{R}^{+}.$$ Taking test functions in the
Schwartz class of smooth functions whose derivatives (up to second order) are rapidly decreasing, we
may assume that
$$f_j, \ f_j'', \ V'f'_j\ \ \ {\rm are\ \ uniformly \ bounded\ \ for\ \ all}\ \
j\geq1.$$ By Ethier and Kurtz \cite{EK} (p.107), to prove the
tightness of $\{L_N(t),\ t\in [0, T],\ N\geq 1\}$,  it is sufficient
to prove that for each $j$ the sequence of continuous real-valued
functions
$$\{\langle L_N(t), f_j\rangle,\ t\in [0, T],\ N\geq1\}$$ is relatively compact.
To this end, note that, by the first part of Theorem \ref{Th1},
there is non-collision and non-explosion for the particles
$\lambda_N^i(t)$ for all $t\in [0, \infty)$. By It\^{o}'s formula,
we have
\begin{eqnarray}\label{LN}
d\langle L_N(t), f\rangle&=&\frac{1}{N}\sqrt{\frac{2}{\beta
N}}\sum\limits_{i=1}^Nf'(\lambda_N^i(t))dW_t^i+ \left\langle L_N(t),
\left(\frac{2}{\beta}-1\right)\frac{1}{2N}f''-\frac{1}{2}V'f'\right\rangle
dt\nonumber\\
& &\hskip2cm
+\frac{1}{2}\int\int_{\mathbb{R}^2}\frac{f'(x)-f'(y)}{x-y}L_N(t,dx)L_N(t,dy)dt.
\end{eqnarray}This yields
\begin{eqnarray}
\langle L_N(t), f_j\rangle&=&\langle L_N(0),
f_j\rangle+\frac{1}{2}\int_0^t\int_{\mathbb{R}}\int_{\mathbb{R}}\frac{f_j'(x)-f_j'(y)}{x-y}L_N(s,dx)L_N(s,dy)ds\nonumber\\
& &-\frac{1}{2}\int_0^t\langle L_N(s), V'f_j'\rangle
ds+\int_0^t\left\langle L_N(s),
\left(\frac{2}{\beta}-1\right)\frac{1}{2N}f_j''\right\rangle ds+M_N^{f_j}(t)\nonumber\\
&=&I_1(N)+I_2(N)+I_3(N)+I_4(N)+M_N^{f_j}(t),\label{E3.A}
\end{eqnarray}
where $$M_N^{f_j}(t)=\frac{1}{N}\sqrt{\frac{2}{\beta
N}}\int_0^t\sum\limits_{i=1}^Nf_j'(\lambda_N^i(s))dW_s^i.$$ Note
that, as $L_N(0)$ is weakly convergent, $I_1(N)$ is convergent. By
the assumption that $f_j$ and $f_j''$ are uniformly bounded (hence
$f_j'$ are uniformly bounded) , we can easily show that
$\{M_N^{f_j}(t), t\in [0, T]\}$ and $I_4(N)$ converge  to zero.
Moreover, by the assumption that $V'f_j'$ and $f_j''$ are uniformly
bounded, the Arzela-Ascoli theorem implies that $I_2(N)$ and
$I_3(N)$ are relatively compact in $C([0, T], \mathbb{R})$. Thus the sequence
$\{(L_N(t))_{t\geq0}:\ N\geq1\}$ is tight in $C([0, T],
\mathbb{R})$.  Tightness also follows for $j=0$ if we have
$$\langle L_N(0),f_0\rangle\rightarrow\ \ {\rm finite\ \ limit\ \ as}\ \ N\rightarrow\infty.$$
So let us suppose that the initial distribution $L_N(0)$ have the
property $\langle L_N(0),f_0\rangle\leq K$ for some $K,$ for all
$N.$ For given $\mu_0,$ we could always find $L_N(0)$ and $f_0$ to
satisfy this and the other conditions, and this gives the tightness
for $j=0$ also.

Finally, we identify the limit process of any weakly convergent
subsequence of $\{L_N(t)\}$. Assuming that $\{L_{N_j}(t), t\in [0,
T]\}$ is a weakly convergent subsequence in $C([0, T],
\mathscr{P}(\mathbb{R}))$. Then,  for all $f\in C_b^2(\mathbb{R})$,
the It\^o's formula $(\ref{E3.A})$ and the above argument show that
$\langle \mu_t,f\rangle=\lim\limits_{j\rightarrow \infty}\langle
L_{N_j}(t),f\rangle$ satisfies the following equation
\begin{eqnarray*}
\int_{\mathbb{R}} f(x)\mu_t(dx)&=&\int_{\mathbb{R}}
f(x)\mu_0(dx)+\frac{1}{2}\int_0^t\int\int_{\mathbb{R}^2}\frac{\partial_xf(x)-\partial_yf(y)}{x-y}\mu_s(dx)\mu_s(dy)ds\\
& &\ \ \ \ \ \ \ -\frac{1}{2}\int_0^t\int_{\mathbb{R}}
V'(x)f'(x)\mu_s(dx)ds.
\end{eqnarray*}
This proves that $\mu_t$ is a weak solution to the McKean-Vlasov
equation $(\ref{DBM7})$.  The proof of Theorem \ref{Th1} is
completed.\hfill $\square$

\section{McKean-Vlasov equation: gradient flow and uniqueness}

To characterize the McKean-Vlasov limit $\mu_t$, we need only to use
the test function $f(x)=(z-x)^{-1}$, where
$z\in\mathbb{C}\backslash\mathbb{R}$, instead of using all test
functions $f\in C_b^2(\mathbb{R})$ in the McKean-Vlasov equation
$(\ref{DBM7})$. Let
\begin{eqnarray*}
G_t(z)=\int_{\mathbb{R}} {\mu_t(dx)\over z-x}
\end{eqnarray*}
be the Stieltjes transform of $\mu_t$. Then $G_t(z)$ satisfies the
following equation
\begin{eqnarray}
{\partial\over \partial t}G_t(z)=-G_t(z) {\partial\over \partial
z}G_t(z)-\frac{1}{2}\int_{\mathbb{R}} {V'(x)\over
(z-x)^2}\mu_t(dx).\label{CHSV-1}
\end{eqnarray}
In particular, in the case $V(x)=\theta x^2$, since
\begin{eqnarray*}
-\int_{\mathbb{R}} {x\over (z-x)^2}\mu_t(dx)=z{\partial\over
\partial z}G_t(z)+G_t(z),
\end{eqnarray*}
the Stieltjes transform of $\mu_t$ satisfies the complex Burgers
equation
\begin{eqnarray}
{\partial\over \partial t}G_t(z)=\left(-G_t(z)+\theta z\right)
{\partial\over
\partial z}G_t(z)+\theta G_t(z).\label{CHSV-2}
\end{eqnarray}

In \cite{Ch, RS93}, Chan and Rogers-Shi proved that the complex
Burgers equation $(\ref{CHSV-2})$ (equivalently, the McKean-Vlasov
equation with potential $V(x)=\theta x^2$) has a unique solution,
and $\lim\limits_{t\rightarrow \infty}G_t(z)$ exists and coincides
with the Stieltjes transform of the Wigner semi-circle law
$\mu_{SC}$. This yields a dynamic proof of the Wigner's theorem,
i.e., $L_N(\infty)$ weakly converges to $\mu_{SC}$.

However, for non quadratic potential $V$, $\int_{\mathbb{R}}
{V'(x)\over (z-x)^2}\mu_t(dx)$ in $(\ref{CHSV-1})$ cannot be
expressed in terms of $G_t(z)$ and its derivatives with respect to
$z$.  Thus, one cannot derive an analogue of the complex Burgers
equation $(\ref{CHSV-2})$  for non quadratic potential $V$, and we
need to find a new approach to prove the uniqueness of the weak
solutions of the Mckean-Vlasov equation for general potential $V$.
In this section, we use the theory of gradient flow on the
Wasserstein space $\mathscr{P}_2(\mathbb{R})$ and the optimal
transportation theory to study this problem.

\subsection{Proof of Theorem \ref{GFV}}

By Theorem \ref{Th1}, we have proved the existence of weak solution
to the McKean-Vlasov equation $(\ref{DBM7})$. Assuming that the weak
solution $\mu_t$ of the McKean-Vlasov equation $(\ref{DBM7})$ is
absolutely continuous with respect to the Lebesgue measure $dx$, we
derive the existence of the weak solution of the nonlinear
Fokker-Planck equation $(\ref{NFK1})$. Thus, to prove the law of
large numbers for $L_N(t)$, we need only to show the uniqueness of
the nonlinear Fokker-Planck equation $(\ref{NFK1})$. Note that,
letting
\begin{eqnarray*}
W(x)=-2\log |x|, \ \ \ \ x\neq 0,
\end{eqnarray*}
then the nonlinear Fokker-Planck equation $(\ref{NFK1})$ can be
rewritten as follows
\begin{eqnarray}
\partial_t \rho=\nabla\cdot (\rho \nabla (V+W*\rho)).\label{MV-0}
\end{eqnarray}

To study the uniqueness and the longtime behavior of the nonlinear
Fokker-Planck equation $(\ref{NFK1})$ (i.e., $(\ref{MV-0})$),  we
first recall Otto's infinite dimensional Riemannian structure on the
Wasserstein space $\mathscr{P}_2(\mathbb{R}^d)$. Fix $fdx\in
\mathscr{P}_2(\mathbb{R}^d)$, the tangent space of
$\mathscr{P}_2(\mathbb{R}^d)$ at $fdx$ is given by
$$
T_{fdx}\mathscr{P}_2(\mathbb{R}^d)=\{sdx: s\in W^{1,
2}(\mathbb{R}^d, \mathbb{R}), \ \ \int_{\mathbb{R}}sdx=0\}.$$ By
\cite{Ot}, for all $s_idx\in T_{fdx}\mathscr{P}_2(\mathbb{R}^d)$,
$i=1, 2$, there exist a unique $p_i\in W^{1,2}(\mathbb{R}^d,
\mathbb{R}^d)$, $i=1, 2$, such that
$$s_i=-\nabla.(f \nabla p_i)$$
In view of this, Otto's infinite dimensional Riemannian metric on
$T_{fdx}\mathscr{P}_2(\mathbb{R}^d)$ is defined by
\begin{eqnarray*}
g_{fdx}(s_1, s_2)=\int_{\mathbb{R}^d}  \langle \nabla p_1, \nabla
p_2\rangle fdx.
\end{eqnarray*}

Next we recall some results and ideas that we borrow from
\cite{CMV}, in which Carrillo, McCann and Villani studied the
following type McKean-Vlasov evolution equation of the granular
media
\begin{eqnarray}
\partial_t \rho=\nabla\cdot (\rho \nabla (\log \rho+V+W*\rho)). \label{MV}
\end{eqnarray}
They proved that the McKean-Vlasov evolution equation can be
realized as a gradient flow of a free energy functional on the
infinite Wasserstein space. More precisely, they proved

\begin{theorem} (Carrillo-McCann-Villani\cite{CMV})\label{th0} Let $V, W$ be
nice functions on $\mathbb{R}^d$, and
\begin{eqnarray}
F(f)=\int_{\mathbb{R}^d} \rho\log \rho dv+\int_{\mathbb{R}^d} \rho V
dv+{1\over 2}\int_{\mathbb{R}^d}\int_{\mathbb{R}^d}
W(x-y)\rho(x)\rho(y)dxdy.\label{F}
\end{eqnarray}
Then the McKean-Vlasov equation $(\ref{MV})$ is the gradient flow of
$F$ with respect to Otto's infinite dimensional  Riemannian metric
on $\mathscr{P}_2(\mathbb{R}^d)$.
\end{theorem}
Moreover, based on Otto's infinite dimensional geometric calculation
on the Wasserstein space, Carrillo, McCann and Villani \cite{CMV}
proved the following entropy dissipation formula

\begin{theorem}(Carrillo-McCann-Villani\cite{CMV}) \label{th1} Denote $\xi:=\nabla(\log \rho+V+W*\rho)$. Then
\begin{eqnarray}
{d\over dt} F(\rho_t)&=&-\int_{\mathbb{R}^n} \rho |\xi|^2dv, \label{ED1}\\
{d^2\over dt^2} F(\rho_t)&=&2\int_{\mathbb{R}^n}\rho {\rm Tr}(D\xi)^T(D\xi)dx+2\int_{\mathbb{R}^n} \langle D^2V\cdot \xi, \xi\rangle\rho dx\nonumber\\
& &+\int_{\mathbb{R}^{2n}} \langle D^2W(x-y)\cdot [\xi(x)-\xi(y)],
[\xi(x)-\xi(y)]\rangle d\rho(x)d\rho(y).   \label{ED2}
\end{eqnarray}
\end{theorem}

\medskip

Inspired by the earlier works due to Biane \cite{Bian03} and
Biane-Speicher \cite{BS01}, and Carrillo-McCann-Villani \cite{CMV},
we can prove the following results, which play a crucial r\^ole in
the proof of the main results of this paper.\footnote{In
\cite{BS01}, Biane and Speicher gave a heuristic proof of the fact
that the probability density of the large $N$-limit of $L_N(t)$
satisfies the McKean-Vlasov equation $(\ref{NFK1})$ (called the free Fokker-Planck equation in \cite{BS01}). Theorem
\ref{Th1} says that $\mu_t$ satisfies the McKean-Vlasov equation
$(\ref{DBM7})$ and integration by parts shows that $\rho_t$
satisfies $(\ref{NFK1})$. Combining this with
Carrillo-McCann-Villani's result in Theorem \ref{th0}, we obtained
Theorem \ref{MVGF} in August 2012.}

\begin{theorem}\label{MVGF}  For all $V: \mathbb{R}\rightarrow [0,
\infty)$ being a $C^2$ function satisfies the condition
$(\ref{Cond1})$, the nonlinear Fokker-Planck equation
$(\ref{NFK1})$, i.e.,
\begin{eqnarray*}
{\partial \rho_t\over \partial t}=-{\partial \over \partial x}(\rho_t({\rm
H}\rho_t-\frac{1}{2}V'))
\end{eqnarray*}
is indeed the gradient flow of  $\Sigma_V$ on the Wasserstein space
$\mathscr{P}_2(\mathbb{R})$.
\end{theorem}
{\it Proof}. Taking $W(x)=2\log|x|^{-1}$, and noting that $\nabla (W*\rho)=H\rho$, Theorem
\ref{th1} follows from Theorem \ref{th0}.   \hfill
$\square$

\begin{theorem} \label{th2} Under the notation of Theorem
\ref{GFV}, we have
\begin{eqnarray}
\label{entro-diss-3}
{d\over dt}\Sigma_V(\mu_t|\mu_V)&=&-2\int_{\mathbb{R}}\left[V'(x)-2{\rm H}\rho_t(x)\right]^2\rho_t(x)dx,\\
{d^2\over dt^2}\Sigma_V(\mu_t|\mu_V)&=&2\int_{\mathbb{R}}  V^{''}(x)|V'(x)-2{\rm H}\rho_t(x)|^2\rho_t(x) dx\nonumber\\
& &+\int\int_{\mathbb{R}^2} {\left[V'(x)-V'(y)-2({\rm
H}\rho_t(x)-{\rm H}\rho_t(y))\right]^2\over (x-y)^2}
\rho_t(x)\rho_t(y)dxdy.
\end{eqnarray}
\end{theorem}
{\it Proof}. By analogue of the proof of Theorem \ref{th1} in
\cite{CMV}, and observing that for $W(x)=-2\log |x|$, we have
$\xi:=\nabla(V+W*\rho) =V'-2{\rm H}\rho$, we can prove Theorem
\ref{th2}. \hfill $\square$

\medskip

\subsection{Proof of Theorem \ref{Th2}}

The proof follows the same argument as used in \cite{Ot, OV, CMV}.
We use the fact that the nonlinear Fokker-Planck equation
$(\ref{NFK1})$ is the gradient flow of the Voiculescu entropy
$\Sigma_V$ on the Wasserstein space $\mathscr{P}_2(\mathbb{R})$, and
that $\Sigma_V$ is $K$-convex along the geodesic displacement
between two probability measures in $\mathscr{P}_2(\mathbb{R})$.

\begin{theorem}\label{K}\footnote{After we proved Theorem \ref{K} in August 2012, we noticed later from Villani's book \cite{Vi2} that Blower \cite{Blo} has proved the $K$-convexity of the Voiculescu entropy.}\label{K-convex} Assuming that $V\in C^2(\mathbb{R},
\mathbb{R}^+)$ and there exists a constant $K\in \mathbb{R}$ such
that $V''\geq K$. Then
$$
{\rm Hess}_{\mathscr{P}_2(\mathbb{R})} \Sigma_V(\mu) \geq K.
$$
\end{theorem}
{\it Proof}. Let $(\rho_s, v_s)$ be the solution to the following continuity equation and the Hamilton-Jacobi equation
\begin{eqnarray*}
\partial_s\rho+\nabla\cdot(\rho v)&=&0,\\
\partial_s (\rho v)+\nabla\cdot(\rho v\otimes v)&=&0.
\end{eqnarray*}
Let $\mu_s=\rho_sdx$.
By analogue of the calculus of the Hessian of the free energy in \cite{CMV}, we can prove that
\begin{eqnarray*}
{d^2\over ds^2}\Sigma_V(\mu_s)=\int_{\mathbb{R}}  V^{''}(x)|v_s(x)|^2\rho_s(x) dx+{1\over 2}\int\int_{\mathbb{R}^2} {|v_s(x)-v_s(y)|^2\over (x-y)^2} \rho_s(x)\rho_s(y)dxdy.
\end{eqnarray*}Thus, under the assumption $V''\geq K$, we have
\begin{eqnarray*}
{\rm Hess}_{\mathscr{P}_2(\mathbb{R})}
\Sigma_V(\mu)(v, v)=\left. {d^2\over ds^2}\Sigma_V(\mu_s)\right|_{s=0}\geq K\int_{\mathbb{R}}  |v(x)|^2\rho_s(x) dx=K\|v\|^2.
\end{eqnarray*} \hfill $\square$

\begin{proposition}\label{MVI}
Let  $\mu(0)(dx)=\rho(0)dx$ and  $\mu(1)(dx)=\rho(1)dx$ be two probability
measures with compact support on $\mathbb{R}$, let
$\mu(s)(dx)=\rho(s)dx$ be the unique geodesic in the Wasserstein space
$\mathscr{P}_2(\mathbb{R})$ linking $\mu(0)$ and $\mu(1)$. Then
\begin{eqnarray}
\left.\left\langle \frac{d\rho(s)}{ds}, {\rm grad}_W \Sigma_V
(\rho(s))\right\rangle\right|_{s= 1} - \left.\left\langle
\frac{d\rho(s)}{ds}, {\rm grad}_W \Sigma_V(\rho(s))
\right\rangle\right|_{s= 0} \geq K W_{2}^{2}(\rho(0),
\rho(1)),\label{DDD1}
\end{eqnarray}
where ${\rm grad}_W\Sigma_V$ denotes the gradient of $\Sigma_V$ on
the Wasserstein space $\mathscr{P}_2(\mathbb{R})$ equipped with
Otto's infinite dimensional Riemannian metric.
\end{proposition}
{\it Proof}.  By the assumptions, we have
\begin{eqnarray*}
{d^2\over ds^2}\Sigma_V(\mu(s)) ={\rm
Hess}_{\mathscr{P}_2(\mathbb{R})}
\Sigma_V(\rho(s))\left(\frac{\partial \rho(s)}{\partial s},
\frac{\partial \rho(s)}{\partial s}\right) \geq
K\left\|\frac{\partial \rho(s)}{\partial
s}\right\|^{2}_{\mathscr{P}_2(\mathbb{R})}.
\end{eqnarray*}
By the mean value theorem, for some $\sigma^*\in (0, 1)$,
\begin{eqnarray*}
\Sigma_V(\rho(1)) - \Sigma_V(\rho(0)) &=& \left.{d\over ds}\right|_{s=0}\Sigma_V(\rho(s)) + {1\over 2}\left.{d^2\over ds^2}\right|
_{s=\sigma^*}\Sigma_V(\rho(\sigma))\\
&\geq& \left.\left\langle \frac{d\rho(s)}{ds}, {\rm grad}_W
\Sigma_V(\rho(s)) \right\rangle\right|_{s= 0} +
{K \over 2}\int^{1}_{0}\left\|\frac{\partial \rho(s)}{\partial s}\right\|^{2}_{\mathscr{P}_2(\mathbb{R})}d\sigma \\
&=&\left.\left \langle \frac{d\rho(s)}{ds}, {\rm
grad}_W\Sigma_V(\rho(s)) \right\rangle\right|_{s= 0}
 + {K \over 2} W_{2}^{2}(\rho(0), \rho(1)).
\end{eqnarray*}

Similarly,
\begin{eqnarray*}
\Sigma_V(\rho(0)) - \Sigma_V(\rho(1)) &\geq& -\left.\left\langle
\frac{d\rho(s)}{ds}, {\rm grad}_W \Sigma_V(\rho(s))
\right\rangle\right|_ {s= 1} + {K\over 2} W_{2}^{2}(\rho(0),
\rho(1)).
\end{eqnarray*}
Summing the two inequalities together, we obtain  $(\ref{DDD1})$.
\hfill $\square$

\medskip
We are ready to give the proof of Theorem \ref{Th2} as follows.\\
{\bf Proof of Theorem \ref{Th2}}. Let $\rho_{t}(s, x)dx:[0, 1]
\rightarrow \mathscr{P}_2(\mathbb{R})$ be the unique geodesic
between $\mu_1(t)$ and $\mu_2(t)$. By Otto \cite{Ot}, we have the
following  derivative formula of the Wasserstein distance
\begin{eqnarray*}
{d \over dt}W^{2}_{2}(\mu_1(t), \mu_2(t)) &=&
-2\int_{\mathbb{R}}\left.\left\langle\frac{d\rho_t(s)}{ds}(x),
\xi_t\right\rangle\right|_{s=0}d\mu_1(t)
 + 2\int_{\mathbb{R}}\left.\left\langle \frac{d\rho_t(s)}{ds}(x), \xi_t\right\rangle\right|_{s=1} d\mu_2(t)\\
&=& 2 \left.\left\langle \frac{d\rho_t(s)}{ds}(x), {\rm grad}_W
\Sigma_V(\mu_2(t))\right\rangle\right|_{s=0}-2\left.\left\langle
\frac{d\rho_t(s)}{ds}(x), {\rm grad}_W
\Sigma_V(\mu_1(t))\right\rangle\right|_{s=1}.
\end{eqnarray*}
By Proposition \ref{MVI}, we have
\begin{eqnarray*}
{d \over dt}W^{2}_{2}(\mu_1(t), \mu_2(t))
 \leq -2KW_{2}^{2}(\mu_1(t),
\mu_2(t)).
\end{eqnarray*}
The Gronwall inequality implies
\begin{eqnarray*}
W_{2}(\mu_1(t), \mu_2(t)) \leq e^{-Kt}W_{2}(\mu_1(0), \mu_2(0)).
\end{eqnarray*}
As a consequence, the McKean-Vlasov equation $(\ref{DBM7})$ has a
unique weak solution. This finishes the proof of Theorem \ref{Th2}.
\hfill $\square$

\subsection{Proof of Theorem \ref{Th3}}
{\bf Proof of Theorem \ref{Th3}}. By Theorem \ref{Th1}, the family
$\{L_N(t), t\in [0, T]\}$ is tight with respect to the weak
convergence topology on $\mathscr{P}(\mathbb{R})$, and the limit of
any weakly convergent subsequence of  $\{L_N(t), t\in [0, T]\}$ is a
weak solution of $(\ref{DBM7})$. By the uniqueness of weak solutions
to $(\ref{NFK1})$, we conclude that $L_N(t)$ weakly converges to
$\mu_t$, and hence $\mathbb{E}[L_N(t)]$ weakly converges to $\mu_t$
as $N\rightarrow \infty$.

Taking $f(x)=x^2$ in $(\ref{DBM7})$ and $(\ref{LN})$ respectively,
we can derive that

\begin{eqnarray}
\label{mut}
\frac{d}{dt}\int_{\mathbb{R}}x^2\mu_t(dx)=1-\int_{\mathbb{R}}
xV'(x)\mu_t(dx),\label{DBM18}
\end{eqnarray}
and
\begin{eqnarray}\label{LN8}
d\langle L_N(t), x^2\rangle&=&\frac{2}{N}\sqrt{\frac{2}{\beta
N}}\sum\limits_{i=1}^N\lambda_N^i(t)dW_t^i+ \left\langle L_N(t),
\left(\frac{2}{\beta}-1\right)\frac{1}{N}-xV'\right\rangle dt+1.
\end{eqnarray}
Taking expectation, we have
\begin{eqnarray}\label{LN8a}
{d\over dt}\int_{\mathbb{R}} x^2 \mathbb{E}[L_N(t, dx)]=1+
\left(\frac{2}{\beta}-1\right)\frac{1}{N}-\int_{\mathbb{R}}xV'(x)L_N(t,
dx).
\end{eqnarray}

On the other hand, from the proof of Theorem \ref{Th1}, we have
\begin{eqnarray*}
\int_{\mathbb{R}} x^2 \mathbb{E}[L_N(t)]=\mathbb{E}\left[{1\over N}\sum\limits_{i=1}^N \lambda_N^i(t)^2\right]=\mathbb{E}[R_t]\leq\mathbb{E}[R'_t].
\end{eqnarray*}
Note that
$$dR'_t={2\over N}\sqrt{R'_t\over \beta}dB_t+\left(\frac{1}{\beta
N}+\frac{N-1}{2N}+\frac{1}{2}\gamma+\gamma R'_t\right)dt,$$
which yields
\begin{eqnarray*}
{d\over dt}\mathbb{E}[R'_t]={1\over \beta N}+\frac{N-1}{2N}+\frac{1}{2}\gamma+\gamma \mathbb{E}[R'_t]\leq {3+\gamma\over 2}+\gamma \mathbb{E}[R'_t].
\end{eqnarray*}
The Gronwall inequality implies
\begin{eqnarray*}
\sup\limits_{t\in [0, T]}\sup\limits_{N}\mathbb{E}[R_t]\leq C(\gamma, \mathbb{E}[R_0]) e^{\gamma T}<\infty.
\end{eqnarray*}
That is
\begin{eqnarray*}
\sup\limits_{t\in [0, T]}\sup\limits_{N}\int_{\mathbb{R}}x^2d\mathbb{E}[L_N](x)\leq C(\gamma, \mathbb{E}[R_0]) e^{\gamma T}<\infty.
\end{eqnarray*}
By H\"older inequality, for all $p\in [1, 2)$,
\begin{eqnarray*}
\int_{|x|\geq A}x^pd\mathbb{E}[L_N(t)](x)\leq \left(\int_{\mathbb{R}}x^2d\mathbb{E}[L_N(t)](x)\right)^{p/2}\left(\mathbb{E}[L_N(t)](|X|\geq A)\right)^{(2-p)/2}.
\end{eqnarray*}
By the tightness of $\mathbb{E}[L_N(t)]$, we have
\begin{eqnarray*}
\lim\limits_{A\rightarrow \infty} \sup\limits_{t\in [0, T]}\sup\limits_{N}\int_{|x|\geq A}x^pd\mathbb{E}[L_N(t)](x)=0.
\end{eqnarray*}
By the characterization of the $W_p$-convergence on
$\mathscr{P}(\mathbb{R})$, see \cite{Vi1, Vi2}, for all $p\in [1, 2)$,  we have
\begin{eqnarray*}
\lim\limits_{N\rightarrow \infty}\sup\limits_{0\leq t\leq T}W_p(\mathbb{E}[L_N(t)], \mu_t)=0.
\end{eqnarray*}
When $V(x)={Kx^2\over 2}$, we have
\begin{eqnarray}
\label{mut}
\frac{d}{dt}\int_{\mathbb{R}}x^2\mu_t(dx)=1-K\int_{\mathbb{R}}
x^2\mu_t(dx),\label{DBM18a}
\end{eqnarray}
\begin{eqnarray}\label{LN8b}
d\langle L_N(t), x^2\rangle&=&\frac{2}{N}\sqrt{\frac{2}{\beta
N}}\sum\limits_{i=1}^N\lambda_N^i(t)dW_t^i-K \left\langle L_N(t),
x^2\right\rangle dt+\left(\frac{2}{\beta}-1\right)\frac{1}{N}+1,
\end{eqnarray}
and
\begin{eqnarray}\label{LN8a}
{d\over dt}\int_{\mathbb{R}} x^2\mathbb{E}[L_N(t, dx)]=1+
\left(\frac{2}{\beta}-1\right)\frac{1}{N}-K\int_{\mathbb{R}}x^2\mathbb{E}[L_N(t,
dx)].
\end{eqnarray}
Hence 
\begin{eqnarray*}
\int_{\mathbb{R}} x^2\mathbb{E}[L_N(t, dx)]-\int_{\mathbb{R}} x^2\mu_t(dx)=e^{-Kt}\left[\int_{\mathbb{R}} x^2\mathbb{E}[L_N(0, dx)]-\int_{\mathbb{R}} x^2\mu_0(dx)\right]+{1\over N}\left(\frac{2}{\beta}-1\right){1-e^{-Kt}\over K}.
\end{eqnarray*}

The proof of Theorem \ref{Th3} is completed.  \hfill $\square$

\subsection{\bf Proof of Theorem \ref{THPC}}
 By the conditions in
Theorem~\ref{Th3} and the Theorem
of~Sznitman~and~Tanaka's~\cite{ST}, we know that,
$M_N(0)$~is~$\mu_0$-chaotic. Since $L_N(t)$~weakly converges to the
deterministic measure~$\mu_t$ for every $t\in [0,T]$, and the
systems~${\rm (GDBM)_V}$ are exchangeable systems, then we have this
propagation of chaos by ~Sznitman~and~Tanaka's Theorem~\cite{ST}.
\hfill $\square$

\medskip

\section{Proof of Theorem \ref{Th5}}

{\bf Proof of Theorem \ref{Th5} (i)}.

By Corollary 3.2 in Biane \cite{Bian03}, for any $C^2$-convex $V$,
there exists a unique equilibrium measure $\mu$ (indeed $\mu=\mu_V$)
with a density $\rho$ satisfying the Euler-Lagrange equation
${\rm H}\rho(x)={1\over 2}V'(x)$ for all $x\in {\rm supp}(\mu)$.  Thus,
$\Sigma_V$ has a unique minimizer $\mu_V$. Moreover,
 as $V$ is $C^2$-convex,  Theorem \ref{th2} implies that
$\Sigma_V$ is a geodesically convex on $\mathscr{P}_2(\mathbb{R})$.

By the fact that $\Sigma_V$ is lower semi-continuous and with
respect to the weak convergence topology, see e.g. \cite{AGZ, Gui},
we see that it is also lower semi-continuous with respect to the
Wasserstein topology on $\mathscr{P}(\mathbb{R})$. Moreover, for all
$c\in \mathbb{R}$ the level set $\{\mu: \Sigma(\mu)\leq c\}$ of
$\Sigma_V$ is relatively compact in the weak convergence topology on
$\mathscr{P}(\mathbb{R})$. By the characterization of the
convergence in the Wasserstein space $\mathscr{P}_2(\mathbb{R})$, we
see that for all $c\in \mathbb{R}$ and $R>0$, $\{\mu:
\Sigma(\mu)\leq c\}\cap B(\mu_0, R)$ is relatively compact with
respect to the topology induced by the Wasserstein distance on
$\mathscr{P}_2(\mathbb{R})$, where $B(\mu_0, R)=\{\mu\in
\mathscr{P}_2(\mathbb{R}): W_2(\mu_0, \mu)\leq R\}$. Hence
$\Sigma_V$ is proper on any geodesic balls of
$\mathscr{P}_2(\mathbb{R})$.

By Proposition $4.1$ in Kloekner \cite{Kloe}, we know that
$\mathscr{P}_2(\mathbb{R})$ has vanishing sectional curvature in the
sense of Alexandrov. More precisely, for any $\mu_1, \mu_2, \mu_3\in
\mathscr{P}_2(\mathbb{R})$ and for any Wasserstein geodesic $\gamma:
[0, 1]\rightarrow \mathscr{P}_2(\mathbb{R})$ such that
$\gamma(0)=\mu_1$ and $\gamma(1)=\mu_2$, for all $t\in [0, 1]$, it
holds that
 \begin{eqnarray*}
 W_2^2(\mu_3, \gamma(t))=tW_2^2(\mu_3, \mu_1)+(1-t)W_2^2(\mu_3, \mu_2)-t(1-t)W_2^2(\mu_1, \mu_2).
 \end{eqnarray*}
Therefore,  $\mathscr{P}_2(\mathbb{R})$ is a nonpositively curved
(NPC) space in the sense of Alexandrov (even though
$\mathscr{P}_2(\mathbb{R}^n)$ is an Alexander space with nonnegative
curvature for $n\geq 2$, see e.g. \cite{AGS}).

By Mayer \cite{Mey} and \cite{Kuw}, we can conclude that $W_2(\mu_t,
\mu_V)\rightarrow 0$ holds if we only assume that $V$ is a
$C^2$-convex potential. The proof of Theorem \ref{Th5} (i) is
completed.\\

\medskip

\noindent{\bf Proof of Theorem \ref{Th5} (ii)}. Taking
$\mu_1(t)=\mu_t$ and $\mu_2(t)\equiv \mu_V$ in Theorem \ref{Th2}, we
have
\begin{eqnarray*}
W_{2}^{2}(\mu_t, \mu_V) \leq e^{-2Kt}W_{2}^{2}(\mu(0), \mu_V).
\end{eqnarray*}
By the fact that $\mu_t$ is the gradient flow of $\Sigma_V$ on
$\mathscr{P}_2(\mathbb{R})$ and using the uniform $K$-convexity of
$\Sigma_V$, we can use the same argument as in \cite{Ot} to prove
\begin{eqnarray*}
\Sigma_V(\mu_t|\mu_V)\leq e^{-2Kt}\Sigma_V(\mu_0|\mu_V).
\end{eqnarray*}
Indeed, by Otto's calculus, we have
\begin{eqnarray*}
\frac{d}{dt}\|{\rm
grad}_W\Sigma_V(\mu_t)\|^2_{\mathscr{P}_2(\mathbb{R})} &=&
2\left\langle {\rm grad}_W\|{\rm
grad}_W\Sigma_V(\mu_t)\|^2_{\mathscr{P}_2(\mathbb{R})}, \frac{d\mu_t}{dt}\right\rangle\\
&=& -2{\rm Hess}_{\mathscr{P}_2(\mathbb{R})} \Sigma_V(\mu_t)\left(\frac{d\mu_t}{dt},  \frac{d\mu_t}{dt}\right)\\
&\leq& -2K\left\|\frac{d\mu_t}{dt}\right\|^2_{\mathscr{P}_2(\mathbb{R})}\\
&=& -2K\|{\rm grad}_W
\Sigma_V(\mu_t)\|^2_{\mathscr{P}_2(\mathbb{R})}.
\end{eqnarray*}
Note that ${\rm grad}_W \Sigma_V(\mu_V)=0$. Thus
\begin{eqnarray*}
\frac{d}{dt}\Sigma_V(\mu_t|\mu_V) &=& \left\langle {\rm
grad}_W \Sigma_V(\mu_t), \frac{d\mu_t}{dt}\right\rangle\\
&=& -\|{\rm grad}_W \Sigma_V(\mu_t)\|^2_{\mathscr{P}_2(\mathbb{R})}\\
&=& \int^{\infty}_{t}\frac{d}{ds}\|{\rm grad}_W \Sigma_V(\mu_s)\|^2_{\mathscr{P}_2(\mathbb{R})}ds\\
&\leq& -2K\int^{\infty}_{t}\|{\rm grad}_W \Sigma_V(\mu_s)\|^2_{\mathscr{P}_2(\mathbb{R})}ds\\
&=& 2K\int^{\infty}_{t}\frac{d}{ds}\Sigma_V(\mu_s)ds\\
&=& -2K\Sigma_V(\mu_t|\mu_V),
\end{eqnarray*}
where in the last step we have used the fact
$\Sigma_V(\mu(\infty))=\Sigma_V(\mu_V)=0$. The Gronwall inequality
implies
\begin{eqnarray*}
\Sigma_V(\mu_t|\mu_V)\leq e^{-2Kt}\Sigma_V(\mu_0|\mu_V).
\end{eqnarray*}
The proof of Theorem \ref{Th5} (ii) is completed. \hfill $\square$

\medskip

To prove Theorem \ref{Th5} (iii), we need the following free
logarithmic Sobolev inequality and free Talagrand transportation
cost inequality due to Ledoux and Popescu \cite{Le-Po09}.

\begin{theorem}\label{LP09}  (Ledoux-Popescu \cite{Le-Po09}) Suppose that $V$ is a $C^2$, convex and there exists a constant $r>0$ such that
$$V''(x)\geq K>0, \ \ \ \ |x|\geq r.$$
Then there exists a constant $c=C(K, r)>0$ such that the free
Log-Sobolev inequality holds: for all probability measure $\mu$ with
$I_V(\mu)<\infty$,
\begin{eqnarray*}
\Sigma_V(\mu|\mu_V)\leq {2\over c}{\rm I}_V(\mu).
\end{eqnarray*}
Moreover, the free Talagrand transportation inequality holds: there
exists a constant $C=C(K, r, V)>0$ such that
\begin{eqnarray*}
CW_2^2(\mu, \mu_V)\leq \Sigma_V(\mu|\mu_V).
\end{eqnarray*}
\end{theorem}
{\bf Proof of Theorem \ref{Th5} (iii)}. By Biane and Speicher
\cite{BS01},  we have the following entropy dissipation formula
\begin{eqnarray*}
{\partial \over \partial t}\Sigma_V(\mu_t|\mu_V)=-{1\over 2 }{\rm I}_V(\mu_t).\label{entro-diss-2}
\end{eqnarray*}
By Theorem \ref{LP09}, there exists a constant $C_1>0$ such that the {free LSI}
holds
\begin{eqnarray*}
\Sigma_V(\mu|\mu_V)\leq {2\over C_1}{\rm I}_V(\mu),
\end{eqnarray*}
which yields
\begin{eqnarray*}
{d\over dt} \Sigma_V(\mu_t|\mu_V)\leq -{C_1\over 4}
\Sigma_V(\mu_t|\mu_V).
\end{eqnarray*}
By the Gronwall inequality, we have
\begin{eqnarray*}
\Sigma_V(\mu_t|\mu_V)\leq e^{-C_1t/4}\Sigma_V(\mu_0|\mu_V).
\end{eqnarray*}
By Theorem \ref{LP09} again, there exists a constant $C_2>0$ such that the  free
transportation cost inequality holds
\begin{eqnarray*}
W_2^2(\mu_t, \mu_V)\leq {1\over C_2}\Sigma_V(\mu_t|\mu_V).
\end{eqnarray*}
Therefore
\begin{eqnarray*}
W_2^2(\mu_t, \mu_V)\leq {e^{-C_1t/4}\over C_2}\Sigma_V(\mu_0|\mu_V).
\end{eqnarray*}
This finishes the proof of Theorem \ref{Th5} (iii). \hfill $\square$

\medskip

\begin{remark} By the same argument as used in Otto \cite{Ot} and Otto-Villani
\cite{OV}, we can prove the following HWI inequality: Suppose that
there exists a constant $K\in \mathbb{R}$ such that
$$V''(x)\geq K, \ \ \ \forall x\in \mathbb{R}.$$
Let $\mu_i\in \mathscr{P}_2(\mathbb{R})$, $i=1, 2$. Then for all
$t>0$, the HWI inequality holds
\begin{eqnarray}
\Sigma_V(\mu_1)-\Sigma_V(\mu_2)\leq W_2(\mu_1, \mu_2)\|{\rm
grad}_W\Sigma_V(\mu_1)\|_{\mathscr{P}_2(\mathbb{R})}-{K\over
2}W_2^2(\mu_1, \mu_2).\label{HWI-1}
\end{eqnarray}
In particular, for any solution to the McKean-Vlasov equation
$(\ref{DBM7})$, we have
\begin{eqnarray}
\Sigma_V(\mu_t)\leq W_2(\mu_t, \mu_V)\|{\rm
grad}_W\Sigma_V(\mu_t)\|_{\mathscr{P}_2(\mathbb{R})}-{K\over
2}W_2^2(\mu_t, \mu_V).\label{HWI-2}
\end{eqnarray}
where
\begin{eqnarray*}
\|{\rm
grad}_W\Sigma_V(\rho)\|^2_{\mathscr{P}_2(\mathbb{R})}=\int_{\mathbb{R}}\rho
|V'(x)-2{\rm H}\rho(x)|^2dx.
\end{eqnarray*}
\end{remark}
To save the length of the paper, we leave the proof to the reader.

\medskip

\section{Double-well potentials and some conjectures}

In this section we discuss again the problem of the longtime
convergence of the McKean-Vlasov equation towards to the equilibrium
measure. More precisely, we want to study the question under which
condition on the external potential $V$ the following double limits
are exchangeable. That is,
\begin{eqnarray*}
\lim\limits_{N\rightarrow \infty}\lim\limits_{t\rightarrow \infty}L_N(t)=
\lim\limits_{t\rightarrow \infty}\lim\limits_{N\rightarrow \infty}L_N(t).
\end{eqnarray*} By \cite{Ch, RS93}, see also \cite{AGZ, Gui}, this is the case when $V(x)={x^2\over 2}$.

Theorem \ref{Th5} ensures the longtime convergence of the weak
solution of the McKean-Vlasov equation to the equilibrium measure
$\mu_V$ for $C^2$-convex potentials $V$. In particular, Theorem
\ref{Th5} applies to {$V(x)=a|x|^{p}$} with $a>0$ and $p\geq 2$.
When $V(x)={x^2\over 2}$ and $\beta=1, 2, 4$, this corresponds to
the cases of GUE, GOE and GSE. Moreover,  Theorem \ref{Th5} also
applies to the Kontsevich-Penner model on the Hermitian random
matrices ensemble with external potential (cf. \cite{CM})
$$
V(x)={ax^4\over 12}-{bx^2\over 2}-c\log|x|.
$$
provided that $a>0, c>0$ and $4ac\geq b^2$.

Can we establish the longtime convergence of the McKean-Vlasov
equation in the non-convex case of external potential?  In
\cite{BS01, Bian03}, Biane and Speicher gave a non-convex potential
$V$ to which the longtime convergence of $\mu_t$ fails. Indeed, as
$\mu_t$ satisfies the gradient flow of the Voiculescu free entropy
$\Sigma_V$ on $\mathscr{P}(\mathbb{R})$, $\mu_t$ may converge to a
local minimizer of $\Sigma_V$ which is not necessary the global
minimizer $\mu_V$. In statistical physics, this indicates that there
might be a phase transition for the large $N$-GDBM model with
non-convex potentials.

Let us consider the double-well potential
\begin{eqnarray*}
V(x)={1\over 4}x^4+{c\over 2}x^2, \ \ \ \ x\in \mathbb{R},\label{DW}
\end{eqnarray*}
where $c\in \mathbb{R}$ is a constant. By
\cite{Joh98, BI}, it has been known that the density function of the
equilibrium measure $\mu_V$ can be explicitly given as follows:  \\

$(i)$ When $c<-2$, $\rho(x)={1\over
2\pi}|x|\sqrt{(x^2-a^2)(b^2-x^2)}1_{[a, b]}$,
where
$a^2=-2-c$ and $b^2=2-c$. \\

$(ii)$ When $c=-2$, $\rho(x)={1\over 2\pi}x^2\sqrt{4-x^2}1_{[-2, 2]}$. \\

$(iii)$ When $c>-2$, $\rho(x)={1\over
\pi}(b_2x^2+b_0)\sqrt{a^2-x^2}1_{[-a, a]}$ ,
where $a^2={\sqrt{4c^2+48}-2c\over 3}$, $b_0={c+\sqrt{{c^2\over
4}+3}\over 3}$, and $b_2={1\over 2}$.\\

When $c\in [0, \infty)$, $V$ is $C^2$ convex and $V''(x)\geq 3$ for $|x|\geq 1$. In this case, Theorem \ref{Th5} (ii) implies that $W_2(\mu_t, \mu_V)\rightarrow 0$ with an exponential convergence rate.

When $c\in(-\infty, -2)$, $\mu_V$ has two supports $[-b, -a]$ and
$[a, b]$ which are disjoint. By Section $7.1$ in Biane-Speicher
\cite{BS01}, it is known that $\mu_t$ does not converge to $\mu_V$.
See also Biane \cite{Bian03}. This also indicates that one cannot
simultaneously prove a free version of the Holley-Stroock
logarithmic Sobolev inequality and a free version of the Talagrand
$T_2$-transportation cost inequality under bounded perturbations of
$p_N(dx)=Z_N^{-1}\prod_{i<j}|x_i-x_j|^2\prod_{i=1}^Ne^{-NV(x_i)}dx$.
Otherwise, by analogue of the proof of Theorem \ref{Th5} (ii), we
may prove that $\mu_t$ converges to $\mu_V$ with respect the
$W_2$-Wasserstein distance and hence in the weak convergence
topology on $\mathscr{P}(\mathbb{R})$. See also \cite{Le-Po09, MMS}
for a discussion on non-convex potentials.

In the case $c\in [-2, 0)$, as the global minimizer $\mu_V$ of
$\Sigma_V$ has a unique support, and all stationary point of $\mu_V$
must satisfy the Euler-Lagrange equation ${\rm H}\mu={1\over 2}V'$,
one can see that the Voiculescu free entropy $\Sigma_V$ has a unique
minimizer which is $\mu_V$. As $\mu_t$ is the gradient flow of
$\Sigma_V$ on $\mathscr{P}_2(\mathbb{R})$, and since ${d\over
dt}\Sigma_V(\mu_t)=-2\int_{\mathbb{R}}\left[V'(x)-2{\rm
H}\rho_t(x)\right]^2\rho_t(x)dx$, we see that $\Sigma_V(\mu_t)$ is
strictly decreasing in time $t$ unless $\mu_t$ achieves the
minimizer $\mu_V$. This yields that the limit of $\Sigma_V(\mu_t)$
exists as $t\rightarrow \infty$. If $\{\mu_t\}$ is tight, and
$\lim\limits_{t\rightarrow \infty}\Sigma_V(\mu_t)=\Sigma_V(\mu_V)$,
we can derive that $\mu_t$ weakly converges to $\mu_V$. By lack of
the tightness of $\{\mu_t\}$, the question whether $W_2(\mu_t,
\mu_V)\rightarrow 0$ (or even $\mu_t$ weakly converges to $\mu_V$)
as $t\rightarrow \infty$ for the above double-well potential $V$
remains open.

We would like to raise the following conjectures.

\begin{conjecture} \label{conj1} Consider the double-well potential $V(x)={1\over 4}x^4+{c\over 2}x^2$ with $c\in [-2, 0)$.
 Then $\mu_t$ converges to $\mu_V$ with respect the $W_2$-Wasserstein distance and hence in the weak convergence topology
 on $\mathscr{P}(\mathbb{R})$.
  \end{conjecture}

\begin{conjecture}\label{conj2}
   Suppose that the potential $V$ is a $C^2$ potential function with $V''(x)\geq K_1$ for all $|x|\geq r$
   and $V''(x)\geq -K_2$ for all $|x|\leq r$, where $K_1, K_2, r>0$ are some positive constants.
   Suppose further that  $\Sigma_V$ has a unique minimizer which has a single compact support. Then $\mu_t$ converges  to $\mu_V$ with respect the $W_2$-Wasserstein distance and in the weak convergence topology on $\mathscr{P}(\mathbb{R})$.
\end{conjecture}

Finally, let us mention the following conjecture due to Biane and
Speicher \cite{BS01}.

\begin{conjecture}\label{conj3} Consider the double-well potential given by $V(x)={1\over 2}x^2+{g\over 4}x^4$, where $g$ is a negative constant but very close to zero. Then $\mu_t$ weakly converges to $\mu_V$.
  \end{conjecture}

\begin{flushleft}
\medskip\noindent
Songzi Li, School of Mathematical Science, Fudan University, 220,
Handan Road, Shanghai, 200432, China, and Institut de Math\'ematiques, Universit\'e Paul Sabatier\\
118, route de Narbonne, 31062, Toulouse Cedex 9, France

\medskip

Xiang-Dong Li, Academy of Mathematics and Systems Science, Chinese
Academy of Sciences, 55, Zhongguancun East Road, Beijing, 100190,
China
\medskip

Yong-Xiao Xie, Academy of Mathematics and Systems Science,
Chinese Academy of Sciences, 55, Zhongguancun East Road, Beijing,
100190, China
\end{flushleft}

\end{document}